\documentclass[12pt]{article}
\usepackage{amsmath, amssymb, amsthm}
\usepackage{graphicx}
\usepackage{mathrsfs}

\usepackage{calc,tikz}
\usepackage{float}
\usetikzlibrary{arrows,decorations.markings,positioning,shapes.geometric}

\tikzstyle{vertex}=[circle, draw, inner sep=1pt, minimum size=4pt]

\tikzstyle{ann} = [fill=white,font=\footnotesize,inner sep=1pt]
\tikzstyle{arrow} = [thick,<-->,>=stealth]

\DeclareMathAlphabet{\mathpzc}{OT1}{pzc}{m}{it}

 \newtheorem{thm}{Theorem}[section]
 \newtheorem{cor}[thm]{Corollary}
 \newtheorem{lem}[thm]{Lemma}
 
 \theoremstyle{definition}
 \newtheorem{defn}[thm]{Definition}

 \theoremstyle{remark}
 \newtheorem{rem}[thm]{Remark}

\numberwithin{equation}{section}

 \newcommand{\address}[1]
 { \vspace{-2em}\begin{center}
  \footnotesize{#1}
   \end{center}}

\title{\Large\bf A Common neighborhood signed graph of a signed graph}
\author{Mukti Acharya}
\date{}

\begin{document}
\maketitle
\address{Former Professor and Head, Department of Applied Mathematics, \\ Delhi Technological University, Delhi- 110042\\
mukti1948@gmail.com}  

\begin{abstract}
Let $G$ be a simple connected graph  with the vertex set $V(G)$ and the edge set $E(G)$ such that $|V(G)| =n $ and $|E(G)|=q$. The common  neighborhood graph $con(G)$ of $G$ has the same set of vertices  as that of $G$ and $uv$ is an edge in $con(G)$ if  $u$ and $v$  have at least one common neighbor in $G$.

Let $S=(G,\sigma)$ be a signed graph, where $G$ is called the underlying graph of $S$ and $\sigma: E(G)\rightarrow \{+,-\}$ is called the signature of $G$. In this paper, the concept of  common neighborhood  graph  has been initiated in the realm of signed graphs, that is, we define a common neighborhood signed graph of a signed graph. \\
We characterize signed $(C^s_n,\mu_\sigma)$ and $(K^s_n,\mu_\sigma)$ so  that   
$con(C^s_n,\mu_\sigma)$ and $con(K^s_n,\mu_\sigma)$ are balanced and also find a solution to  $con(C^s_n,\mu_\sigma)$ $\sim$ $(C^s_n,\mu_\sigma)$ and $con(K^s_n,\mu_\sigma)$ $\sim$$(K^s_n,\mu_\sigma)$, where $\sim$ is a switching  equivalence relation.    

.

\end{abstract}
\hrulefill 

\section{Introduction}
The graphs considered in this paper are assumed to be connected and simple unless otherwise mentioned. Let $G$ be a graph with vertex set $V(G)$ and edge set $E(G)$ having $|V(G)| =n$ and $|E(G)|=q$. Let $S=(G,\sigma)$ be a signed graph, where $G$ is called the underlying graph of $S$ and $\sigma: E(G)\rightarrow \{+,-\}$ is called the signature of $G$. The edges that  receive sign $+$  are called positive edges, denoted by $E^{+}(S)$, and the edges that receive sign  $-$  are called negative edges, denoted by $E^{-}(S)$. The set of edges $E(S)=E^{+}(S)\cup E^{-}(S)$. In signed graphs, the positive edges are drawn as solid lines, and the negative edges are drawn as dashed lines.  
 
\indent A signed graph in which all the edges are positive is called an all-positive signed graph, and a signed graph in which all the edges are negative is called an all-negative signed graph. An all-positive (all-negative) signed graph is called homogeneous; otherwise, it is called heterogeneous. For more information on signed graphs, we refer to \cite{TZ}. \\ 
 \indent The  marking of the vertices of $S$ defined as 
          $\mu_\sigma(v) =  \prod \sigma(vu)$ for $vu$ in $E(S)$ 
is called a canonical marking.The vertices receiving the sign $+$ are called positive vertices, and the vertices receiving the sign $-$ are called negative vertices.\\ 
\indent By $d(v)$ we mean the total degree of a vertex $v\in V(S)$. By $d^{+}(v)$ and $d^{-}(v)$, we denote the number of positive and negative edges incident with $v\in V(S)$ . Thus, $d(v) = d^{+}(v) + d^{-}(v)$, where $d^{+}(v)$ is called the positive degree of $v$ and  $d^{-}(v)$ is called the negative degree of $v$. By $N_G(v)$, we denote the set of  vertices of $G$ adjacent to a vertex $v$ and will be called the set of vertices in the neighborhood of $v$. For any two vertices $x$ and $y$, by $A = N_G(x) \cap N_G(y)$ we denote the common neighborhoods of $x$ and $y$. 
Similarly, we can define the common neighborhoods of $x$ and $y$ in $S$.\\
\indent A function  $\mu:V(S)\rightarrow\{+, -\}$ is called the marking of the signed graph $S$. Switching a signed graph $S$ with respect to a given marking $\mu$, we obtain another signed graph $S^{*}$ from $(S,\mu)$ by changing the sign of its edges $uv$ for which $\mu(u)\neq\mu(v)$, where by $(S,\mu)$ we mean a signed graph $S$ along with its marking $\mu$. A signed graph $S$ is said to switch to another signed graph $S^{'}$, written as $S\sim S^{'}$, whenever there exists a marking $\mu$ such that $(S,\mu)\sim S^{'}$. It is obvious that $\sim$ is an equivalence relation in  the class of all signed graphs. Hence, if $S\sim S^{'}$, we will mean that $S$ and $S^{'}$ are switching equivalent. For terminology used in the paper and not defined here, we refer to \cite{FH2}.
Now we give a few definitions that we need for our work
\begin{defn}\cite{MA}
Let $S$ be a signed graph. By a positive (negative) section of $S$ we mean a maximal connected induced subgraph in $S$ consisting of only positive (negative) edges of $S$.
\end{defn}

\begin{defn}\cite{FH1}
A signed graph $S$ is said to be balanced if every cycle  has an even number of negative edges.
\end{defn}
\begin{defn}\cite{FH2}
 The line graph $L(G)$ of $G=(V,E)$ has $E(G)$ as its vertex set and two vertices are joined by an edge in $L(G)$ if the corresponding edges have a vertex in common in $G$.  
\end{defn}
The following theorem will be found useful for our work.
\begin{thm}\cite{FH1} A signed graph S is balanced  if and only if its set of vertices $V(S)$ can be partitioned into two subsets $V_1(S)$ and $V_2(S)$, one of them may be empty, such that negative edges lie across $V_1(S)$ and $V_2(S)$.
\end{thm}

\section{2-path graph}
\begin{defn}[\cite{ EMR}] The 2-path graph $P_2(G)$ of $G$ has the same set of vertices  as that of G and two vertices $u$ and $v$ are joined by an edge $uv$ in $P_2(G)$ if there is a 2-path between $u$ and $v$ in $G$.
\end{defn}
A graph $G$ for which \\
\indent{}      $P_n(G)\cong G$ ...........................................(1)\\
has been termed an $n-path$ invariant graph in \cite{EM} and  for  $n=2$, we will have an $2-path$ invariant graph, where $\cong$ denotes the isomorphism.  For a detailed study on the notion of n-path graphs, we refer to \cite { A, EMR, EM1,EM,EH, FHK, S}.\\
\indent The following theorem characterizes the structure of 2-path invariant graphs.
\begin{thm}[\cite{EM}] \label{t:thm 2.2 } 
A graph $G$ of order $p$ is a 2-path invariant graph if and only if $G\cong \overline {K_p}$ or $K_p $ with $p\geq3$, or $C_p$, $p=2m+1, m\geq1$.
\end{thm}

\section{2-path signed graph}
\begin{defn}[\cite{GP}] The 2-path signed graph $P_2(S)$ of $S$ has the same set of vertices  as that of S and two vertices $u$ and $v$ are joined by an edge $uv$ in $P_2(S)$, if there is a 2-path between $u$ and $v$, and the edge $uv$ is negative if all 2-paths are all-negative and otherwise the edge $uv$ is positive.
\end{defn}
A signed graph $S$ for which \\
\indent{}      $P_n(S)\sim S$ ..............................................(2)\\
has been termed the  switching  invariant $n-path$  signed graph in \cite{G,GP} and  for  $n=2$, we will have the  switching  invariant $2-path$  signed  graph, where $\sim$ denotes the switching equivalent relation. For a detailed study on the notion of an  n-path signed graphs, we refer to \cite{MA1,G,GP,DD}. 

\begin{defn}[\cite{SAAM,PAL}]  The common  neighborhood graph $con(G)$ of $G$ has the same set of vertices  as that of $G$ and $uv$ is an edge in $con(G)$ if  $u$ and $v$  have at least one common neighbor in $G$.
\end{defn}
\begin{rem}  How are $P_2(G)$ and $con(G)$ of $G$ related ? Both concepts are equivalent.
\end{rem}
\section{ Common neighborhood signed graph}

 We now define the common neighborhood signed graph of a signed graph.  
\begin{defn} Let $S=(G,\sigma)$ be a signed graph and $\mu_\sigma (u)=\prod_{v\in N_S(u)} \sigma(uv)$. The common neighborhood signed graph  $con(S,\mu_\sigma )$ of $(S,\mu_\sigma )$ has the same set of vertices  as that of $(S,\mu_\sigma )$ and $uv$ is an edge in $con(S,\mu_\sigma)$ if  $u$ and $v$  have at least one common neighbor in $(S,\mu_\sigma )$ and $\sigma_{con}(uv)= \prod_{w \in A} \mu_\sigma(w)$.

\indent{}The readers may like to refer to \cite{AA} for the definition of a line signed graph taking into account the  marking  of a common vertex of adjacent edges, and by marking we mean the canonical marking. 
\end{defn}
\begin{rem} How are $P_2(S)$ and $con(S,\mu_\sigma)$ of $S$ related ? The signature of the edges in $P_2(S)$ and $con(S,\mu_\sigma)$ differs.
\end{rem}
Note that when $\sigma$ is all-positive, this recovers a classical common neighborhood  graph.

Instead of using a 2-path signed graph, we will adhere to a  common neighborhood signed graph  throughout the paper.  
\section{ Signed Cycles}
Now we explore the structure of signed cycles  $(C^s_n,\mu_\sigma)$ of order  $n$ and the vertices have canonical markings, so that  $con(C^s_n,\mu_\sigma)$ is  balanced. The following theorem will be useful for our further studies.
\begin{thm} For  $G=C_n$, $n\geq3$, $con(G)$ is connected if $n$ is odd and $con(G)$ is disconnected having exactly two components, each is a cycle of length $\frac{n}{2}$ if $n$ is even.  
\end{thm}

\begin{thm}
For  $(C^s_n,\mu_\sigma)$, $n=2m+1, m\geq1 $,  $con(C^s_n,\mu_\sigma)$ is always balanced.     
\end{thm}
\begin{proof}
By definition, it is clear that  every negative section contributes two negative edges in $con(C^s_n,\mu_\sigma)$. Thus, the number of negative edges in $con(C^s_n,\mu_\sigma)$ is even.Therefore, $con(C^s_n,\mu_\sigma)$ is always balanced.
\end{proof}
\begin{cor} For homogeneous $(C^s_n,\mu_\sigma)$, $n = 2m+1, m\geq1 $, $con(C^s_n,\mu_\sigma)$ is balanced.
\end{cor}
\begin{thm} 
For even $n$, $con(C^s_n,\mu_\sigma)$ is balanced if and only if the number of negative vertices having odd labels is  even in number in $(C^s_n,\mu_\sigma)$.
\end{thm}
\begin{proof}
Let $u_i, i=1,2,...,n$, be the vertices of $(C^s_n,\mu_\sigma)$ for even $n$. According to  Theorem 5.1, $con(C^s_n,\mu_\sigma)$ is disconnected having exactly two components and  each is a cycle of length $\frac{n}{2}$. We assume that $u_1,u_3,u_5...., u_{n-1}$  are the vertices of cycle $C1$ and $u_2,u_4,u_6...., u_n$ are the vertices of cycle $C2$ and  $(C^s_n,\mu_\sigma)$ has an even number of negative vertices. Now, we have the following two cases\\   
Case 1: $C1$ and $C2$ in $(C^s_n,\mu_\sigma)$ have an even number of negative vertices.\\ 
Clearly, each cycle in $con(C^s_n,\mu_\sigma)$ has an even number of negative edges. Hence, $con(C^s_n,\mu_\sigma)$ is balanced.\\
Case 2:  $C1$ and $C2$ in $(C^s_n,\mu_\sigma)$ have an odd number of negative vertices.\\
Clearly, each cycle in $con(C^s_n,\mu_\sigma)$ has an odd number of negative edges. Hence, $con(C^s_n,\mu_\sigma)$ is not balanced. The converse is easy to see. 
Therefore, the result follows.
\end{proof}

\begin{lem} If $n\cong2(mod)4$ and each  negative and positive section in $(C^s_n,\mu_\sigma)$ is of length $1$ then $con(C^s_n,\mu_\sigma)$ is not balanced.   
\end{lem}
\begin{proof} Clearly, $d^-(u_i)=1, u_i=1,2,...,n$ in $(C^s_n,\mu_\sigma)$ , that is, in the neighborhood of any two vertices in $(C^s_n,\mu_\sigma)$, we have a  negative vertex. Thus, both the cycles $C1$ and $C2$ in $con(C^s_n,\mu_\sigma)$ are negative homogeneous of odd length. Therefore, $con(C^s_n,\mu_\sigma)$ is not balanced. 
\end{proof}
\section{ Complete Signed Graphs}

Now, we explore the structure of $(K^s_n,\mu_\sigma)$ of order $n$ when its vertices have canonical marking so that  $con(K^s_n,\mu_\sigma)$ is  balanced.
\begin{thm} For $(K^s_n,\mu_\sigma)$, $con(K^s_n,\mu_\sigma)$ is always balanced.
\end{thm}
\begin{proof} To prove the result, we  partition $V(K^s_n,\mu_\sigma)$ into two subsets, namely $V_1(K^s_n,\mu_\sigma)$ having positive vertices  and $V_2(K^s_n,\mu_\sigma)$ having  negative vertices. Now, we show that $\bigl \langle V_1(K^s_n,\mu_\sigma)\bigr \rangle$ and $\bigl\langle V_2(K^s_n,\mu_\sigma)\bigr\rangle$ are positive homogeneous in $con(K^s_n,\mu_\sigma)$, where  $\bigl\langle V(S)\bigr\rangle$ denotes the induced subgraph on the set of vertices of $S$.
To prove this, we consider the following cases:\\
Case I: Assume that $uv \in \bigl\langle V_1(K^s_n,\mu_\sigma)\bigr\rangle$ and $W$ is the set of  common neighbors of $u$ and $v$  $\in (K^s_n,\mu_\sigma)$ and $|W|= n-2$. Clearly, $\sigma_{con}(uv)=+$ for $uv\in con(K^s_n,\mu_\sigma)$ as in the common neighborhood of $u,v\in (K^s_n,\mu_\sigma)$  we have  an even number of negative vertices. Therefore, $\bigl\langle V_1(K^s_n,\mu_\sigma)\bigr\rangle$ is positive homogeneous. Note that $uv \in \bigl\langle E( V_1(K^s_n,\mu_\sigma))\bigr\rangle$ can be a positive or a negative edge.\\
Case II: Using the same arguments as in Case I, we conclude that $\bigl\langle V_2(K^s_n,\mu_\sigma)\bigr\rangle$ is positive homogeneous.\\ 
To prove that $con(K^s_n,\mu_\sigma)$ is balanced, we must  show that the edges between $V_1$ and $V_2$ in $(K^s_n,\mu_\sigma)$ are negative in $con(K^s_n,\mu_\sigma)$. Let $uv\in E(K^s_n,\mu_\sigma)$, where $u\in V_1$ has $\mu_\sigma(u)=+$ and $v\in V_2$ has $\mu_\sigma(v)=-$. Assuming  that $uv \in E(K^s_n,\mu_\sigma)$ is a positive edge  and in the neighborhood of $u$ and $v$ in $uv\in(K^s_n,\mu_\sigma)$, we have an odd number of vertices that are negative and, hence, $uv\in(K^s_n,\mu_\sigma)$ becomes a negative edge in  $ con(K^s_n,\mu_\sigma)$. Using the same arguments, we conclude that the negative edge 
$uv\in(K^s_n,\mu_\sigma)$ remains negative in $ con(K^s_n,\mu_\sigma)$. Therefore, $ con(K^s_n,\mu_\sigma)$ is balanced.\\
Now we show that one of the subsets $V_1$ and $V_2$ in $(K^s_n,\mu_\sigma)$ may be empty. Assume that the subset $V_1$ is empty. Then, clearly $|V_2| = n$ and the vertices in $V_2$ are negative. Therefore, every $uv\in E(K^s_n,\mu_\sigma)$ becomes a positive edge in $ con(K^s_n,\mu_\sigma)$. Therefore, $ con(K^s_n,\mu_\sigma)$ is balanced. Hence, the result follows. 
\end{proof}

\begin{cor}
If $(K_n^s,\mu_\sigma)$ is homogeneous,  then $con(K_n^s,\mu_\sigma)$ is  balanced.
\end{cor}
\begin{cor}
Let $(K_n^s,\mu_\sigma)$ be of odd order. If $M$ is negative homogeneous non perfect matching and the remaining edges of  $(K_n^s,\mu_\sigma)$ are positive then $con(K_n^s,\mu_\sigma)$ is balanced.
\end{cor}
\section{ Switching invariant common neighborhood signed graphs}
A signed graph $S$ for which \\
\indent     $ con(S)\sim S$ ............................(3)\\
will be  termed the  switching  invariant  common neighborhood  signed  graph, where $\sim$ denotes the switching equivalent relation; that is,  we have to characterize $S$ so that $con(S)$ is switching equivalent to $S$\\
\indent{}In view of Theorem \ref{t:thm 2.2 }, we see that if $S$ is a solution to (3) then $S$ is a totally disconnected graph, or  $S\in (K^s_n,\mu_\sigma)$ for an integer $n\geq3$, or $S\in (C^s_{2m+1},\mu_\sigma)$ for an integer $m\geq1$. Therefore, in order to completely find the structure of switching invariant common neighborhood signed graphs, it is enough to search for solutions for (3) in the sets $(C^s_{2m+1},\mu_\sigma)$ for a positive integer $m\geq1$, and $(K^s_n,\mu_\sigma)$ for an  integer $n\geq3$.\\

The following theorem is needed for our work.
\begin{thm}\cite{Z} \label{thm: t 6.1}
For a given graph $G$ any two signed graphs $S_1$ and $S_2$ whose underlying graph is $G$ are switching equivalent if and only if they are cycle-isomorphic.
\end{thm}
The signed graphs $S_1$ and $S_2$ are said to be cycle-isomorphic if there exists an isomorphism $f$ from $G_1$ onto $G_2$, where $G_1$ and $G_2$ are the underlying graphs of $S_1$ and $S_2$, respectively, such that the sign of every cycle $Z$ in  $S_1$ equals the sign of $f(Z)$ in $S_2$, where the sign of any subgraph is defined as the product of the signs of its edges. The following theorem determines the solutions to (3) that can be found in 
$(C^s_{2m+1},\mu_\sigma)$ for some positive integer $m\geq1$.
\begin{thm} $Z \in (C^s_{2m+1},\mu_\sigma)$, for any integer $m\geq1$, is a solution to (3) if and only if $Z$ is balanced.
\end{thm}
\begin{proof} We have already shown that for $Z \in (C^s_{2m+1},\mu_\sigma)$, for any integer $m\geq1$, $con(Z)$ is balanced. From Theorem \ref{thm: t 6.1}, it follows that $Z$ must be balanced.The converse is easy to see. Thus, the result follows.  
\end{proof}
\begin{thm}  $(K^s_n,\mu_\sigma)$ for some positive integer $n\geq3$ is a solution to (3) if and only if $(K^s_n,\mu_\sigma)$ is balanced.
\end{thm} 
\begin{proof} We have already shown that $con(K^s_n,\mu_\sigma)$ for a positive integer $n\geq3$ is balanced. Therefore, $(K^s_n,\mu_\sigma)$ will be a solution to (3) if $(K^s_n,\mu_\sigma)$ is balanced.The converse is easy to see. Thus, the proof follows.
\end {proof}

\section{Conclusion and further scope}
In this paper, a study  on common neighborhood signed graphs $con(S)$ of $S$ has been initiated. We have characterized $(C^s_n,\mu_\sigma) $ for  $n=2m+1, m\geq1$,  and  $(K^s_n,\mu_\sigma)$ for $n\geq3$, so that $con(C^s_n,\mu_\sigma)$ and $con( K^s_n,\mu_\sigma)$ are balanced.  Solutions to $ con(S)\sim S$ have  also been found for  $(C^s_n,\mu_\sigma)$  and  $(K^s_n,\mu_\sigma)$. The properties of common neighborhood  signed graphs for other classes of signed graphs and its relation to derived signed graphs have to be investigated. We pose a few problems to nurture the study of common neighborhood graphs in the realm of signed graphs.\\
Problem 1: Characterize $S$ so that $con(S,\mu_\sigma)$ is balanced.\\
Problem 2: Characterize $G$ so that $L(G)\cong con(G)$.\\
Problem 3: Characterize $S$ so  that $L(S)\sim con(S,\mu_\sigma)$.\\  
Problem 4: Characterize $S$ so that $P_2(S)$ and $con(S,\mu_\sigma)$  are negation of each other.\\ 
\section{Acknowledgment}
 The author thanks Professor Thomas Zaslavsky for his continuous thought-provoking comments for refining ideas and  providing encouragement.


\begin{thebibliography}{10}\label{bibliography}
\bibitem{A}
B.D. Acharya, Contributions to the theories of hypergraphs, graphoids and graphs, Doctoral dissertation, 1975, IIT, Bombay.
\bibitem{MA1}
Mukti Acharya, Switching invariant three-path signed graphs,  In: Optimization, Design of Experiments and Graph Theory (Eds.: M.N. Gopalan and G.A. Patwardhan), Indian Institute of Technology, Bombay, 1988, 342-345
\bibitem{AA}
 B.D. Acharya and M. Acharya, Dot-line signed graphs, Annals of Pure and Applied Mathematics, 10(1) (2015), 21-27.

\bibitem{MA}
M. Acharya, C-Cordial labeling in signed graphs, Electronic Notes in Discrete Math., 63,(2017), 11--22

\bibitem{EMR}
F. Escalante, L. Montajano and T. Rojano, A charcterization of n-path graphs and of graphs having nth root, J. Combin. Theory Ser. B 16(3), 1974, 282-289. 
\bibitem{EM1}
F. Escalante and L. Montajano, Extremal   n-path connected graphs, Commun. Tecnicas Ser. B CIMAS, National University Maxico, 4(60), 1973. 

\bibitem{EM}
F. Escalante and L. Montajano, Trees and n-path invariant graphs, Graph Theory Newsletter, 3(3),1974. 
\bibitem{EH}
G. Exoo, F. Harary, Step graphs, J. of Combinatorics, Information and System Sciences, 5 (1980), pp 52 - 53.

\bibitem{G}
M K Gill, Contribution to some topics in graph theory and its applications, Ph.D. Thesis, Indian Institute of Technology, Bombay, 1983.
\bibitem{GP}
M.K. Gill and G.A. Patwardhan, Switching invariant two-path signed graphs, Discrete Mathematics 61(1986), 189-196.  
\bibitem{FH1}
F. Harary, On the notion of balance of a signed graph, Michigan Mathematical Journal, 2(1955), 113-146. 

\bibitem{FH2}

F. Harary, Graph Theory, Addison-Wesley, 1972.

\bibitem{FHK}
F. Harary, C. Hoede and D. Kadlecek, Graph value functions related to step graphs, 7(3)(1982), 231-245.  

\bibitem{SAAM}
Samaneh Hossain-Zadeh, Ali Iranmanesh, Asma Hamzeh and Mohammad Ali Hosseinzadeh, On the common neighborhood graphs, Electronic Notes in Discrete Mathematics, 45(2014), 51-56. 

\bibitem{JAG} 
J. A. Gallian, \emph{A dynamic survey of graph labeling}, The Electronic Journal of Combinatorics \textbf{1} (2018),  \# DS6. 
Available: http://www.combinatorics.org 
\bibitem{PAL}
Jeepamol J. Palathingal, V. Anu and S. Aparna Lakshman, A note on common neighborhood graphs,  AIP Conference Proc. 2773, 060001(2023).

\bibitem{S}
S.K. Simic, Graph equations for line graphs and nth distance graphs, Publ. Instt. Math.
(Beograd) 33(47)(1983), 203-216.

\bibitem{DD}
Deepa Sinha and Deepakshi Sharma, Characterization of 2-Path Product Signed Graphs with
Its Properties, Computational Intelligence and Neuroscience
Volume 2017, Article ID 1235715, 8 pages, https://doi.org/10.1155/2017/1235715

\bibitem{Z}
T. Zaslavsky, Signed graphs, Discrete Appl. Math. 4(1982),47-74 ( Erratum:ibid, 5 (1983),248). 

\bibitem{TZ}
T. Zaslavsky, \emph{A mathematical bibliography of signed and gain graphs and allied areas}, Electron. J. Combin. Dynamic Surveys \#DS8, 1998 et seq.
\end{thebibliography}
\end{document}